\documentclass[12pt]{article}

\usepackage{amsmath,amssymb,amsthm}
\usepackage{amsfonts}
\usepackage[mathscr]{eucal}
\begin{document}
\newtheorem{thm}{Theorem}
\numberwithin{thm}{section}
\newtheorem{lemma}[thm]{Lemma}
\newtheorem{remark}{Remark}
\newtheorem{corr}[thm]{Corollary}
\newtheorem{proposition}{Proposition}
\newtheorem{theorem}{Theorem}[section]
\newtheorem{deff}[thm]{Definition}
\newtheorem{case}[thm]{Case}
\newtheorem{prop}[thm]{Proposition}
\numberwithin{equation}{section}
\numberwithin{remark}{section}
\numberwithin{proposition}{section}
\newtheorem{corollary}{Corollary}[section]
\newtheorem{others}{Theorem}
\newtheorem{conjecture}{Conjecture}\newtheorem{definition}{Definition}[section]
\newtheorem{cl}{Claim}
\newtheorem{cor}{Corollary}
\newcommand{\ds}{\displaystyle}

\newcommand{\stk}[2]{\stackrel{#1}{#2}}
\newcommand{\dwn}[1]{{\scriptstyle #1}\downarrow}
\newcommand{\upa}[1]{{\scriptstyle #1}\uparrow}
\newcommand{\nea}[1]{{\scriptstyle #1}\nearrow}
\newcommand{\sea}[1]{\searrow {\scriptstyle #1}}
\newcommand{\csti}[3]{(#1+1) (#2)^{1/ (#1+1)} (#1)^{- #1
 / (#1+1)} (#3)^{ #1 / (#1 +1)}}
\newcommand{\RR}[1]{\mathbb{#1}}
\thispagestyle{empty}
\begin{titlepage}
\title{\bf Lifetime asymptotics of iterated Brownian motion in $\RR{R}^{n}$}
\author{Erkan Nane\thanks{ Supported in part by NSF Grant \# 9700585-DMS }\\
Department of Mathematics\\
Purdue University\\
West Lafayette, IN 47906 \\
enane@math.purdue.edu}
\maketitle
\begin{abstract}
\noindent {\it Let  $\tau _{D}(Z) $ be the first exit time of
iterated Brownian motion from a domain $D \subset \RR{R}^{n}$
started at $z\in D$ and let $P_{z}[\tau _{D}(Z)  >t]$ be its
distribution.   In this paper
 we establish the exact asymptotics of $P_{z}[\tau _{D}(Z)  >t]$
 over bounded domains as an improvement of the results in
 \cite{deblassie, nane2}, for $z\in D$
\begin{eqnarray}
 \lim_{t\to\infty}
t^{-1/2}\exp(\frac{3}{2}\pi^{2/3}\lambda_{D}^{2/3}t^{1/3})
P_{z}[\tau_{D}(Z)>t]= C(z),\nonumber
\end{eqnarray}
where $C(z)=(\lambda_{D}2^{7/2})/\sqrt{3 \pi}\left(
\psi(z)\int_{D}\psi(y)dy\right) ^{2}$. Here $\lambda_{D}$ is the
first eigenvalue of the Dirichlet Laplacian $\frac{1}{2}\Delta$ in
$D$, and $\psi $ is the eigenfunction corresponding to
$\lambda_{D}$ .

 We also study  lifetime asymptotics of Brownian-time Brownian
 motion (BTBM),
$Z^{1}_{t}=z+X(|Y(t)|)$, where $X_{t}$ and $Y_{t}$ are independent
one-dimensional Brownian motions.}

\end{abstract}

\end{titlepage}
\section{Introduction and statement of main results}

Iterated Brownian motion (IBM) have attracted the interest of
several authors \cite{allouba2,  allouba1, bandeb,
burdzy1, burdzy2, bukh, csaki, deblassie, eisenbaumshi, klewis,
 nane, nane2, nane5, xiao}.  Several other iterated processes including
 Brownian-time Brownian motion (BTBM) have also been studied \cite{allouba2, allouba1,
 koslew, nane3, nane4}. One of the main differences between these iterated processes and Brownian motion is
that they  are not Markov processes.  However, these processes
have many properties similar to that of Brownian motion (see
\cite{allouba1,bandeb, deblassie, nane}, and references therein).

To define iterated  Brownian motion $Z_{t}$
started at $z \in \RR{R}$,
 let $X_{t}^{+}$, $X_{t}^{-}$ and $Y_{t}$
 be three  independent
one-dimensional Brownian motions, all started at $0$. Two-sided
Brownian motion is defined by
\[ X_{t}=\left\{ \begin{array}{ll}
X_{t}^{+}, &t\geq 0\\
X_{(-t)}^{-}, &t<0.
\end{array}
\right. \] Then iterated Brownian motion started at $z \in
\RR{R}$ is
\[ Z_{t}=z+X(Y_{t}),\ \ \    t\geq 0.\]

In $\RR{R}^{n}$, one requires $X^{\pm}$ to be independent
$n-$dimensional Brownian motions. This is the version of the
iterated Brownian motion due to Burdzy, see \cite{burdzy1}.

We next define another closely related process, the so called
 Brownian-time Brownian motion. Let $X_{t}$ and $Y_{t}$
 be two  independent
one-dimensional Brownian motions, all started at $0$.
Brownian-time Brownian motion started at $z\in \RR{R}$ is
$$ Z_{t}^{1}=z+X(|Y_{t}|)\ \ \    t\geq 0.$$
In $\RR{R}^{n}$ one requires $X$ to be an $n-$dimensional Brownian
motion.

Let $\tau_{D}$ be the first exit time of Brownian motion from a
domain $D\subset \RR{R}^{n} $. The large time behavior of
$P_{z}[\tau_{D}>t]$ has been studied for  several types of domains,
 including general cones
\cite{bansmits, deblassie1}, parabola-shaped domains \cite{bds,
lshi}, twisted domains \cite{DSmits}, unbounded convex domains \cite{li} and bounded domains
\cite{blab}. Our aim in this article is to do the same for the first
exit time of IBM over bounded domains in $\RR{R}^{n}$, and for the
first exit time of BTBM over several domains in  $\RR{R}^{n}$. See
Ba\~{n}uelos and DeBlassie \cite{bandeb}, Li \cite{li}, Lifshits
and Shi \cite{lshi} and Nane \cite{nane} for a survey of results
obtained for Brownian motion and iterated Brownian motion in these
domains.

For many bounded domains  $D\subset \RR{R}^{n} $ the asymptotics
of $P_{z}[\tau_{D}>t]$ is well-known (See \cite{blab} for a more
precise statement of this.) For $z\in D$,
\begin{equation}\label{exp}
\lim_{t\to \infty} e^{\lambda_{D} t}P_{z}[\tau_{D}>t]=
\psi(z)\int_D \psi(y) dy,
\end{equation}
where $\lambda_{D}$ is the first eigenvalue of $\frac{1}{2}
\Delta$ with Dirichlet boundary conditions  and $\psi$ is its
corresponding eigenfunction.

DeBlassie \cite{deblassie} proved that for iterated Brownian
motion in bounded domains; for $z\in D$,
\begin{equation}\label{boundeddomain0}
\lim_{t\to \infty}t^{-1/3}\log P_{z}[\tau_{D}(Z)>t]=\ -\frac{3}{2}
\pi^{2/3} \lambda_{D}^{2/3}.
\end{equation}

The limits (\ref{exp}) and (\ref{boundeddomain0}) are very
different in that the latter involves taking the logarithm which
may kill many unwanted terms in the exponential. It is then
natural to ask if it is possible to obtain an analogue of
(\ref{exp}) for IBM. That is, to remove the log in
(\ref{boundeddomain0}). In \cite{nane2}, we improved the limit in
(\ref{boundeddomain0}) as follows; for $z\in D$,
\begin{eqnarray}
 2C(z)  & \leq &
\liminf_{t\to\infty} t^{-1/2}
\exp(\frac{3}{2}\pi^{2/3}\lambda_{D}^{2/3}t^{1/3})P_{z}[\tau_{D}(Z)>t]\nonumber\\
& \leq  & \limsup_{t\to\infty}
t^{-1/2}\exp(\frac{3}{2}\pi^{2/3}\lambda_{D}^{2/3}t^{1/3})
P_{z}[\tau_{D}(Z)>t]\leq \pi C(z),\nonumber
\end{eqnarray}
where $ C(z)=\lambda_{D}\sqrt{2\pi /3}\left(
\psi(z)\int_{D}\psi(y)dy\right) ^{2}. $

In this paper we prove the following theorem which improves both
limits above.

\begin{theorem}\label{boundeddomain1}
Let $D\subset \RR{R}^{n}$ be a domain for which (\ref{exp})
holds pointwise and let $\lambda_D$ and $\psi$ be as above. Then
for $z\in D$,
\begin{eqnarray}
 \lim_{t\to\infty}
t^{-1/2}\exp(\frac{3}{2}\pi^{2/3}\lambda_{D}^{2/3}t^{1/3})
P_{z}[\tau_{D}(Z)>t]= \frac{(\lambda_{D}2^{7/2})}{\sqrt{3
\pi}}\left( \psi(z)\int_{D}\psi(y)dy\right) ^{2}.\nonumber
\end{eqnarray}
\end{theorem}

\begin{remark}
Observe that $2 \lambda_{D}\sqrt{2\pi /3}\leq
(\lambda_{D}2^{7/2})/\sqrt{3 \pi} \leq \pi \lambda_{D}\sqrt{2\pi
/3}  $, so Theorem \ref{boundeddomain1} is in agreement with the
results obtained previously in \cite{deblassie, nane2}.
\end{remark}

In \cite{DSmits}, DeBlassie and Smits studied the tail behavior of
the first exit time of Brownian motion in twisted domains in the
plane. Let $D\subset \RR{R}^{2}$ be a domain whose boundary
consists of three curves (in polar coordinates)
\begin{eqnarray*}
\ & C_{1}:&   \theta =f_{1}(r), \ \  r\geq r_{1}\\
\ & C_{2}:&  \theta = f_{2}(r), \ \ r\geq r_{1} \\
\ & C_{3}:&   r=r_{1},  \ \ \ \ \ \  f_{2}(r)\leq \theta  \leq
f_{1}(r)
\end{eqnarray*}
where $f_{1}$ and $f_{2}$ are smooth and the curves $C_{1}$ and
$C_{2}$ do not cross:
$$
0<f_{1}(r)-f_{2}(r)<\pi, \ \ r\geq r_{1}.
$$
DeBlassie and Smits call $D$ a twisted domain if there
 are constants $r_{0}>0$, $\gamma >0$ and $p\in (0,1]$ and a smooth function $f(r)$ such
that the curves $f_{1}(r)$ and $f_{2}(r)$, $r\geq r_{0},$ are
obtained from $f(r)$  by moving out $\pm \gamma r^{p}$ units along
the normal to the curve $\theta=f(r)$ at the point whose polar
coordinates are $(r,f(r))$. They call $\gamma r^{p}$ the growth
radius and $\theta =f(r)$ the generating curve. DeBlassie and
Smits \cite[Theorem 1.1]{DSmits} have the following tail behavior
of the first exit time of Brownian motion in twisted domains
$D\subset \RR{R}^{2}$ with growth radius $\gamma r^{p}$, $\gamma
>0, $ $0<p<1$
\begin{equation}\label{brownianlimit2}
\lim_{t\to \infty}t^{-(\frac{1-p}{1+p})}\log
P_{z}[\tau_{D}>t]=-l_{1}=-\left[ \frac{\pi^{2p-1}}{\gamma
2^{2p}(1-p)^{2p}} \right]^{\frac{2}{p+1}} C_{p}
\end{equation}
where
$$
C_{p}=(1+p)\left[
\frac{\pi^{2+p}}{8^{p}p^{2p}(1-p)^{1-p}}\frac{\Gamma ^{2p} \left(
\frac{1-p}{2p}\right)}{\Gamma ^{2p}\left( \frac{1}{2p}\right)}
\right]^{\frac{1}{p+1}}.$$
 For these domains, Nane \cite{nane} obtained the following; for all $z\in D$,
 $$ \lim_{t\to\infty} t^{-(\frac{1-p}{3+p})}\log
P_{z}[\tau_{D}(Z)>t]= -(\frac{3+p }{2+ 2 p } )
(\frac{1+p}{1-p})^{(\frac{1-p}{3+p})} \pi ^{(\frac{2-2p}{3+p } )}
l_{1}^{  (\frac{2+2p }{3+p }) } ,
$$
where $l_{1}$ is the limit given by (\ref{brownianlimit2}).

We obtained in \cite{nane2}, the following for BTBM in twisted
domains; for $z\in D$,
$$ \lim_{t\to\infty} t^{-(\frac{1-p}{p+3})}\log
P_{z}[\tau_{D}(Z^{1})>t]= -2^{(\frac{2p -2}{3+p})}(\frac{3+p }{2+
2 p } ) (\frac{1+p}{1-p})^{(\frac{1-p}{3+p})} \pi
^{(\frac{2-2p}{3+p } )} l_{1}^{  (\frac{2+2p }{3+p }) } ,
$$
where $l_{1}$ is the limit given by the limit given by
(\ref{brownianlimit2}).

DeBlassie and Smits \cite{DSmits} also obtained similar results
for $p=1$. Let $D\subset \RR{R}^{2}$ be a twisted domain with
growth
 radius $\gamma r$, $\gamma >0$. Then for $z\in D$,
\begin{equation}\label{p=1}
\lim_{t\to\infty}[\log t]^{-1}\log P_{z}[\tau_{D}>t]=-C(\gamma
)=-\pi \left[ 4 \arccos\frac{1 }{\sqrt{1+\gamma
^{2}}}\right]^{-1}.
\end{equation}

We obtain the following lifetime asymptotics of BTBM in twisted
domains.
\begin{theorem}\label{twisted1}
Let $D\subset \RR{R}^{2}$ be a twisted domain with growth
 radius $\gamma r$, $\gamma >0$. Then for $z\in D$,

 $$
\lim_{t\to\infty}[\log t]^{-1}\log
P_{z}[\tau_{D}(Z^{1})>t]=-C(\gamma )/2,
 $$
where $C(\gamma )$ as in (\ref{p=1}).
\end{theorem}

Using Theorem 1.3. from \cite{nane2}, which says that for all
$z\in D$ and all $t>0$,
$$
P_{z}[\tau_{D}(Z)>t]\ \leq \ 2 P_{z}[\tau_{D}(Z^{1})>t],
$$
 we obtain the
following for IBM in twisted domains.
\begin{corollary}
Let $D\subset \RR{R}^{2}$ be a twisted domain with growth
 radius $\gamma r$, $\gamma >0$. Then for $z\in D$,
$$
\limsup_{t\to\infty}[\log t]^{-1}\log P_{z}[\tau_{D}(Z)>t]\leq
-C(\gamma )/2,
 $$
where $C(\gamma )$ as in (\ref{p=1}).
\end{corollary}

In \cite{li}, using Gaussian techniques, Li studied lifetime
asymptotics of Brownian motion in domains of the following form
 $$P_{f}=\{(x,y)\in \RR{R}^{n+1}:y>f(x), x\in
\RR{R}^{n}\}$$ for $f(x)=\exp(|x|^{p}), p>0$. Li established that for
$z\in P_{f}$,
\begin{equation}
\lim_{t\to\infty}t^{-1}(\log t)^{2/p}\log
P_{z}[\tau_{P_{f}}>t]=-j^{2}_{\nu},
\end{equation}
where $\nu =(n-2)/2$ and $j_{\nu}$ is the smallest positive zero
of the Bessel function $J_{v}.$

We obtain the following theorem in these domains
\begin{theorem}\label{theorem3}
Let $P_{f}$ be as above with $f(x)=\exp(|x|^{p}), p>0$. Then for
$z\in P_{f}$,
$$
\lim_{t\to\infty}t^{-1/3}(\log t)^{4/3p}\log
P_{z}[\tau_{P_{f}}(Z^{1})>t]=-C(p),
$$
where
$$
C(p)=(3/2)^{(4+3p)/3p}2^{1/3}(j_{\nu}^{2}2^{2/p})^{2/3}(\pi
^{2}/8)^{1/3}.
$$
\end{theorem}
Using Theorem 1.3. from \cite{nane2}, we obtain the following for
IBM in these domains.
\begin{corollary}
Let $P_{f}$ be as above with $f(x)=\exp(|x|^{p}), p>0$. Then for
$z\in P_{f}$,
$$
\limsup_{t\to\infty}t^{-1/3}(\log t)^{4/3p}\log
P_{z}[\tau_{P_{f}}(Z)>t]\leq -C(p).
$$
\end{corollary}

 For $f(x)=h(|x|)$ and $h^{-1}(x)=Ax^{\alpha}(\log
x)^{\beta}$, $x>1$. Li obtained the following: let $\epsilon>0$.
For $t$ large, $z\in P_{f}$
\begin{equation}
-(1+\epsilon)C_{\alpha,\beta, 1} \leq
t^{-\frac{(1-\alpha)}{(1+\alpha)}}(\log t)^{\frac{2\beta
}{(1+\alpha)}}\log P_{z}[\tau_{P_{f}}>t] \leq -(1-\epsilon)
C_{\alpha,\beta, 2},.
\end{equation}
where
$$
C_{\alpha,\beta, 1}=2^{-1}(1-\alpha)^{-1}(\alpha ^{-\alpha}(1+\alpha)^{2\beta +2}A^{-2}j_{\nu}^{2})1/(1+\alpha)
$$
 and
 $$
 C_{\alpha,\beta, 2}=(1+\alpha)(2\alpha )^{-\alpha /(1+\alpha)}(2^{-1}(1+\alpha))^{2\beta /(1+\alpha)}) C^{1/(1+\alpha)}
 $$
 where $C=(1-\alpha)^{-1}2^{2\beta -1}A^{-2}j_{\nu}^{2}. $

We have the following for BTBM in these domains.
\begin{theorem}\label{not-sharp}
For $0<\alpha <1$  and $\beta \in \RR{R}$,
\begin{eqnarray}
-C(1)&\leq &  \liminf_{t\to\infty} t^{-(1-\alpha)/(3+ \alpha )}
(\log t)^{4\beta (1+\alpha)/(3 +\alpha) }\log
P_{z}[\tau_{P_{f}}(Z^{1})>t]\nonumber\\
& \leq  & \limsup_{t\to\infty} t^{-(1-\alpha)/(3\alpha +1)} (\log
t)^{4\beta (1+\alpha)/(3+\alpha) }\log
P_{z}[\tau_{P_{f}}(Z^{1})>t]
 \nonumber\\
&\leq & -C(2)\nonumber
\end{eqnarray}
where $$C(1)=\left( \frac{3+\alpha}{2(1+\alpha)}
\right)^{\frac{3+\alpha+4\beta}{3+\alpha}}\left(
\frac{1-\alpha}{(3+\alpha)}\right)
(\pi^{2}/8)^{\frac{1-\alpha}{(3+\alpha)}}(C_{\alpha ,\beta ,
1})^{\frac{2(1+\alpha)}{(3+\alpha)}}2^{\frac{4\beta}{(3+\alpha)}},
$$ and
$$C(2)= \left( \frac{3+\alpha}{2(1+\alpha)} \right)^{\frac{3+\alpha+4\beta}{3+\alpha}}\left( \frac{1-\alpha}{(3+\alpha)}\right)
(\pi^{2}/8)^{\frac{1-\alpha}{(3+\alpha)}}(C_{\alpha ,\beta ,
2})^{\frac{2(1+\alpha)}{(3+\alpha)}}2^{\frac{4\beta}{(3+\alpha)}}.$$
\end{theorem}

Using Theorem 1.3. from \cite{nane2}, we obtain the following for
IBM in these domains.
\begin{corollary}For $0<\alpha <1$  and $\beta \in \RR{R}$.
Let $P_{f}$ be as above with $f(x)=h(|x|) $,
$h^{-1}(x)=Ax^{\alpha}(\log x)^{\beta}$, $x>1$. Then for $z\in
P_{f}$,
$$
\limsup_{t\to\infty} t^{-(1-\alpha)/(3\alpha +1)} (\log t)^{4\beta
(1+\alpha)/(3+\alpha) }\log P_{z}[\tau_{P_{f}}(Z^{1})>t]
 \leq  -C(2)
$$
where $C(2)$ is as above.
\end{corollary}

The paper is organized as follows. In \S 2, we give some
preliminary lemmas to be used in the proof of main results.
Theorem \ref{boundeddomain1} is proved in \S3. \S4 is devoted to
prove Theorems \ref{twisted1}, \ref{theorem3} and \ref{not-sharp}.
In \S5, we recall several asymptotic results to be used in the
proof of main results from Nane \cite{nane2}.
\section{Preliminaries}
In this section we state some preliminary facts that will be used
in the proof of main results.

 In what follows  we will write $f
\approx g$ and $f \lesssim g$ to mean that for some positive
$C_{1}$ and $C_{2}$, $C_{1}\leq f/g \leq C_{2}$ and $f \leq C_{1}
g$, respectively. We will also write $f(t) \sim g(t)$, as
$t\rightarrow \infty $,  to mean that $f(t) / g(t) \rightarrow 1$,
as $t\rightarrow \infty $.

 The main fact is the following Tauberian theorem (\cite[Laplace transform method, 1958, Chapter 4]{debruijn}). Laporte \cite{laporte} also studied this type
of integrals.
Let $h$  and $f$ be continuous functions on $\RR{R}$. Suppose $f$
is non-positive and has a global max at $x_{0}$, $f'(x_{0})=0$, $
f''(x_{0})<0$ and $h(x_{0})\neq 0$  and $\int_{-\infty}^{\infty} h(x)\exp(\lambda f(x)) <\infty$ for all $\lambda>0$. Then as $\lambda \to\infty$,
\begin{equation}\label{laplacemethod}
\int_{0}^{\infty}h(x)\exp(\lambda f(x))dx \sim
h(x_{0})\exp(\lambda f(x_{0}))\sqrt{\frac{2\pi}{\lambda
|f''(x_{0})|}}.
\end{equation}

It can be easily seen from Laplace transform method that as $\lambda \to \infty$,
\begin{equation}\label{generalasymptotic}
\int_{0}^{\infty}\exp (-\lambda (x+x^{-2}))dx\sim \exp(-3 \lambda
2^{-2/3})\sqrt{\frac{2^{4/3}\pi}{3 \lambda}}.
\end{equation}

Similarly, as $t \to \infty$,
\begin{equation}\label{newasymptotic}
\int_{0}^{\infty}\exp (-\frac{at}{u^{2}}-bu)du\ \sim
\sqrt{\frac{\pi}{3}}2^{2/3}a^{1/6}b^{-2/3}t^{1/6} \exp (-3
a^{1/3}b^{2/3}2^{-2/3}t^{1/3}).
\end{equation}

This  follows from equation
(\ref{generalasymptotic}) and after making the change of variables
$u=(atb^{-1})^{1/3}x$.

Finally, we obtain, as $t \to \infty$,
\begin{equation}\label{generalast11}
\int_{0}^{\infty}u\exp (-\frac{at}{u^{2}}-bu)du\ \sim 2\sqrt{\frac{\pi}{
3}}a^{1/2} b^{-1} t^{1/2}\exp (-3 a^{1/3}b^{2/3}2^{-2/3}t^{1/3}).
\end{equation}

Writing the power series for the cosine function we easily see
that as $t\to\infty$,
\begin{eqnarray}
&  & \int_{ 0}^{\infty} x\cos \left(\pi K / x \right) \exp
(-\frac{\pi ^{2} t}{2 x^{2}}-\lambda_{D} x)dx \nonumber\\
 & &\sim
2\sqrt{\frac{\pi}{3}}(\frac{\pi^{2}}{2})^{1/2}\lambda_{D}^{-1}
t^{1/2} \exp
(-\frac{3}{2}\pi^{2/3}\lambda_{D}^{2/3}t^{1/3}).\label{generalast111}
\end{eqnarray}
Next we give a Lemma that will be used in the proof of Theorem \ref{twisted1}.
\begin{lemma}\label{upper-bound}
Let $\xi$ be a positive random variable such that as $x\to 0^{+}$
$$[|\log x|]^{-1}\log P[\xi \leq x]\sim -c/2.$$ Then as $\lambda \to \infty$
$$
[\log \lambda]^{-1}\log E[\exp(-\lambda \xi)] \sim -c/2.
$$

\end{lemma}
\begin{proof}
Let $\epsilon >0$, then by hypothesis there exists $\delta
(\epsilon)
>0$ such that
$$
\delta ^{c(1+\epsilon)/2} \leq P[\xi \leq \delta ]\leq \delta
^{c(1-\epsilon)/2},
$$
for all $\delta < \delta (\epsilon )$. Let $f$ be the density of
$\xi$. Then
\begin{eqnarray}
E[\exp(-\lambda \xi)]&= &\int_{0}^{\delta }e^{-\lambda x} f(x)dx
+\int_{\delta}^{\infty}e^{-\lambda x} f(x)dx\nonumber\\
&\leq & P[\xi\leq \delta] +\exp (-\delta \lambda)\nonumber
\end{eqnarray}
Now we use the fact that $\exp(-x)\leq c_{N}x^{-N}$ for any $N\in
\RR{N}$ and for some $c_{N}>0$.

Hence
\begin{eqnarray}
E[\exp(-\lambda \xi)] &\leq & \delta ^{c(1-\epsilon)/2} +
c_{N}(\delta \lambda)^{-N}.\nonumber
\end{eqnarray}

To minimize this upper bound we require

$$
\delta ^{c(1-\epsilon)/2}=c_{N}(\delta \lambda)^{-N}
$$
which gives

$$
\delta=(c_{N}(\lambda)^{-N})^{1/(N+c(1-\epsilon)/2)}.
$$
Hence for some $D(N)>0$,
$$
E[\exp(-\lambda \xi)]\leq D(N) \lambda
^{-\frac{Nc(1-\epsilon)/2}{(N+c(1-\epsilon)/2)}}
$$
Taking logarithm of both sides, dividing by $\log \lambda$, and
letting $\lambda\to\infty$, we obtain
$$
\limsup_{t\to\infty}[\log \lambda]^{-1}\log E[\exp(-\lambda
\xi)]\leq -\frac{Nc(1-\epsilon)/2}{(N+c(1-\epsilon)/2)}
$$
Now letting $N\to\infty$, and  $\epsilon\to 0$, we get
$$
\limsup_{t\to\infty}[\log \lambda]^{-1}\log E[\exp(-\lambda
\xi)]\leq -c/2.
$$

Lower bound follows from

\begin{eqnarray}
E[\exp(-\lambda \xi)] & \geq & \int_{0}^{\delta }e^{-\lambda x}
f(x)dx\nonumber\\
&\geq & e^{-\delta \lambda}P[\xi \leq \delta] \geq \delta
^{c(1+\epsilon)/2}e^{-\delta \lambda}\nonumber
\end{eqnarray}
and taking $\delta =\lambda ^{-1}$, for large $\lambda $, we get
$$
E[\exp(-\lambda \xi)]\gtrsim \lambda ^{-c(1+\epsilon)/2}.
$$
Taking logarithm of both sides and dividing by $\log \lambda$, and
letting $\epsilon\to 0$, we get
$$
\liminf_{t\to\infty}[\log \lambda ]^{-1}E[\exp(-\lambda \xi)]\geq
-c/2.
$$

\end{proof}

We next state a version of de Bruijn's Tauberian Theorem (Kasahara
\cite[Theorem 3]{kasahara} and Bingham, Goldie and Teugels
\cite[p. 254]{bgt}.)
\begin{theorem}[de Bruijn Tauberian Theorem]\label{debruijn}
Let $\xi$ be a positive random variable. Then, for $\alpha >0$ and
$\beta\in \RR{R}$
$$
\log P[\xi \leq \epsilon] \sim -C\epsilon^{-\alpha} |\log
\epsilon|^{\beta} \ \ \ as \ \ \ \epsilon\to 0^{+}
$$
if and only if
$$
\log E[\exp(-\lambda \xi)]\sim
-(1+\alpha)^{1-\beta/(1+\alpha)}\alpha^{-\alpha
/(1+\alpha)}C^{1/(1+\alpha)}\lambda^{\alpha /(1+\alpha)}(\log
\lambda )^{\beta /(1+\alpha)}
$$
as $\lambda\to\infty.$
\end{theorem}

We give next an application of de Bruijn's Tauberian Theorem.
\begin{lemma}\label{densitylemma}
Let $\xi$ be a positive random variable with density $f(x)=\gamma
e^{-v}V$, $v(x)=Cx^{-1/2}(\log x^{-1/2})^{-2/p} $, and $dv=-Vdx$
 Then as $x\to 0^{+}$
 $$
 P[\xi\leq x]\sim -Cx^{1/2}(|\log x^{1/2}|)^{-2/p}.
 $$
 In this case
 $$
 \log E[\exp(-\lambda \xi)]\sim -(3/2)^{(4+3p)/3p}2^{1/3}(C2^{2/p})^{2/3}\lambda^{1/3}(\log \lambda )^{-4/3p}
 $$
\end{lemma}

\section{Iterated Brownian motion in bounded domains}
If $D\subset \RR{R} ^{n}$  is an open set, write
\[
\tau_{D}^{\pm}(z)=\inf \{ t\geq 0:\ \ X_{t}^{\pm} +z \notin D\},\]
and if $I\subset \RR{R}$ is an open interval, write
\[
\eta _{I}=\eta (I)= \inf \{ t\geq 0 :\ \  Y_{t}\notin I\}.
\]
Recall that $\tau _{D}(Z)$ stands for the first exit time of
iterated Brownian motion from  $D$. As in DeBlassie
\cite[\S3.]{deblassie}, we have by the continuity of the paths for
$Z_{t}=z+X(Y_{t})$, if $f$ is the probability density of $\tau
_{D}^{\pm}(z)$
\begin{equation}\label{translation}
P_{z}[\tau _{D}(Z) > t]=\int_{0}^{\infty} \!\int_{0}^{\infty}
 P_{0}[\eta _{(-u,v)}>t]  f(u) f(v)dvdu.
\end{equation}

\begin{proof}[\textbf{The proof of Theorem \ref{boundeddomain1} }]
The following is well-known
\begin{equation}\label{distr1}
P_{0}[\eta _{(-u,v)} >t]=\frac{4}{\pi
}\sum_{n=0}^{\infty}\frac{1}{2n+1} \exp (-\frac{(2n+1)^{2}\pi
^{2}}{2(u+v)^{2}} t) \sin \frac{(2n+1)\pi u}{u+v},
\end{equation}
(see Feller \cite[pp. 340-342]{feller}).\newline Let $\epsilon
>0$. From Lemma \ref{lemmaA.1}, choose $M>0$ so large that
\begin{equation}\label{papprox1}
 (1-\epsilon)  \frac{4}{\pi} e ^{- \frac{\pi ^{2} t}{2} } \sin \pi x \leq P_{x}[\eta _{(0,1)} >t]
 \leq (1+\epsilon)  \frac{4}{\pi} e ^{- \frac{\pi ^{2} t}{2} } \sin \pi x,
\end{equation}
 for $t\geq M $, uniformly $ x \in (0,1)$.
For a bounded domain with regular boundary it is well-known (see
\cite[page 121-127]{blab}) that there exists an increasing
sequence of eigenvalues, $\lambda_{1}< \lambda_{2} \leq
\lambda_{3}\cdots ,$ and eigenfunctions $\psi_{k}$ corresponding
to $\lambda_{k}$ such that,
\begin{equation}\label{eigenvalue-ex0}
P_{z}[\tau_{D}\leq t]=\sum_{k=1}^{\infty} \exp (-\lambda_{k}
t)\psi_{k}(z)\int_{D}\psi_{k}(y)dy.
\end{equation}
From the arguments in DeBlassie \cite[Lemma A.4]{deblassie}
\begin{equation}\label{eigenvalue-ex}
f(t)=\frac{d}{dt}P_{z}[\tau_{D}\leq t]=\sum_{k=1}^{\infty}
\lambda_{k} \exp (-\lambda_{k} t)\psi_{k}(z)\int_{D}\psi_{k}(y)dy.
\end{equation}
 Finally choose $K>0$ so large that
 $$A(z)(1-\epsilon)\exp (-\lambda_{D} u)\leq f(u)\leq A(z)(1+\epsilon)\exp (-\lambda_{D} u) $$
for all $u\geq K$, where
$$A(z)=\lambda_{1}\psi_{1}(z)\int_{D}\psi_{1}(y)dy=\lambda_{D}\psi (z)\int_{D}\psi (y)dy.$$
We further assume that $t$ is so large that $K<\frac{1}{2}
\sqrt{t/M}$.
 Define $A$ for $ K >0 $ and $ M>0$ as
\[
A=\left\{(u,v):\ \ K\leq u\leq \frac{1}{2} \sqrt{\frac{t}{M}} ,\
 u \leq v\leq \sqrt{\frac{t}{M}}-u
\right\}.
\]

 By equation (\ref{papprox1}) and from equation (3.10) in \cite{deblassie},
  \begin{eqnarray}
& & P_{z}[\tau _{D}(Z) > t]=2\int_{0}^{\infty} \!\int_{u}^{\infty}
P_{\frac{u}{u+v}} [\eta
_{(0,1)}>\frac{t}{(u+v^{2})}] f(u)f(v)dvdu\nonumber\\
 &\geq & C^{1} \int_{K}^{\frac{1}{2} \sqrt{t/ M}}
\int_{u}^{\sqrt{t/ M}-u} \sin \left(\frac{\pi u }{
(u+v)}\right)\exp (-\frac{\pi ^{2} t}{2 (u+v)^{2}})
 \exp (-\lambda_{D} (u+v)) dvdu, \nonumber
\end{eqnarray}
where $C^{1}=C^{1}(z)=2(4/\pi)A(z)^{2}(1-\epsilon)^{3}$. Changing
the variables $x=u+v, z=u$ the integral is
\[
= C^{1}\int_{K}^{\frac{1}{2} \sqrt{t/ M}} \int_{2z}^{\sqrt{t/
M}}\sin \left(\frac{\pi z }{ x}\right) \exp(-\frac{\pi ^{2} t}{2
x^{2}}) \exp (-\lambda_{D} x) dxdz,
\]
and reversing the order of integration
\begin{eqnarray}
\ &\ & =C^{1}\int_{2K}^{\sqrt{t/ M}} \int_{K}^{\frac{1}{2} x}\
\sin \left(\frac{\pi z }{ x}\right)\exp(-\frac{\pi ^{2} t}{2
x^{2}}) \exp (-\lambda_{D}
x)dzdx \nonumber \\
 \ &\ & = C^{1}/\pi \int_{2K}^{\sqrt{t/ M}}  x\cos \left(\frac{\pi K }{ x}\right)\exp(-\frac{\pi ^{2} t}{2 x^{2}})
\exp (-\lambda_{D} x)dx \nonumber
\end{eqnarray}

From equation (\ref{generalast111}) as $t\to \infty$,
\begin{eqnarray}
& &\int_{ 0}^{\infty} x\cos \left(\frac{\pi K }{ x}\right) \exp
(-\frac{\pi ^{2} t}{2 x^{2}}) \exp (-\lambda_{D} x)dx \nonumber\\
&\sim &
2\sqrt{\frac{\pi}{3}}(\frac{\pi^{2}}{2})^{1/2}\lambda_{D}^{-1}
t^{1/2} \exp
(-\frac{3}{2}\pi^{2/3}\lambda_{D}^{2/3}t^{1/3}).\label{exactast}
\end{eqnarray}
Now for some $c_{1}>0$,
\begin{eqnarray}
\ &\ &
\int_{ 0}^{ K/ \delta}  x\exp (-\frac{\pi ^{2}
t}{2x^{2}}-\lambda_{D} x)
 dx
\nonumber  \\
\ &\ & \leq  e^{- \pi^2 \delta ^{2}t /2K^{2}}\int_{ 0}^{ K/
\delta}x\exp (-\lambda_{D} x)dx
\lesssim  e^{-c_{1} t},\label{exactast1}
\end{eqnarray}
and
\begin{eqnarray}
\ & \ & \int_{ \sqrt{t/M}}^{\infty} x \exp (-\frac{\pi ^{2} t}{2
x^{2}}) \exp (-\lambda_{D} x)dx
 \leq  \int_{ \sqrt{t/M}}^{\infty}
x\exp (-\lambda_{D} x)dx  \nonumber \\
  \ &\ & = ( \sqrt{t/M} \lambda_{D}^{-1}+\lambda_{D}^{-2})\exp (-\lambda_{D}\sqrt{t/M} ). \label{exactast2}
\end{eqnarray}
Now from equations (\ref{exactast})-(\ref{exactast2}) we get
\begin{equation}\label{exactlowerbound1}
\liminf_{t\to\infty}  t^{-1/2} \exp
(\frac{3}{2}\pi^{2/3}\lambda_{D}^{2/3}t^{1/3}) P_{z}[\tau
_{D}(Z)>t]\geq (C^{1}/\pi)
2\sqrt{\frac{\pi}{3}}(\frac{\pi^{2}}{2})^{1/2}\lambda_{D}^{-1}.
\end{equation}
For the upper bound for $P[\tau _{D}(Z)>t]$ from equation (3.10)
in \cite{deblassie},
\begin{equation}
P_{z}[\tau _{D}(Z) > t]=2\int_{0}^{\infty} \!\int_{u}^{\infty}
P_{\frac{u}{u+v}} [\eta _{(0,1)}>\frac{t}{(u+v^{2})}]
f(u)f(v)dvdu\ .
\end{equation}
 We define the following sets that make up the domain of integration,
\begin{eqnarray}
 & A_{1}& =\{(u,v):v\geq u\geq 0, \  u+v\geq \sqrt{t/M}
 \},\nonumber\\
 & A_{2} & =\{ (u,v): \ u\geq 0,\ v\geq K, \ u\leq v,  \  u+v\leq \sqrt{t/M}
 \},\nonumber\\
& A_{3} & =\{ (u,v):\  0 \leq u\leq v \leq K\}.\nonumber
\end{eqnarray}
Over the set $A_{1}$ we have for some $c>0$,
\begin{eqnarray}
\ & \ & \int\!\int_{A_{1}} P_{\frac{u}{u+v}} [\eta _{(0,1)}>\frac{t}{(u+v)^{2}}] f(u)f(v)dvdu\nonumber\\
\ &\ & \leq \int\!\int_{A_{1}}  f(u)f(v)dvdu\leq
\exp(-c\sqrt{t/M}).\label{upperA1}
\end{eqnarray}
The equation (\ref{upperA1}) follows from the distribution of
$\tau_{D}$ from Lemma 2.1 in \cite{nane}.

Since on $A_{3}$, $t/(u+v)^{2}\geq M$,
\begin{eqnarray}
\ & \ & \int\!\int_{A_{3}} P_{\frac{u}{u+v}} [\eta _{(0,1)}>\frac{t}{(u+v)^{2}}] f(u)f(v)dvdu\nonumber\\
\ &\ & \leq \int_{0}^{K}\!\int_{0}^{K} \exp (-\frac{\pi ^{2} t}{2 (u+v)^{2}}) f(u)f(v)dvdu.\nonumber\\
\ &\ & \leq
\exp(-\frac{\pi^{2}t}{8K^{2}})\int_{0}^{K}\!\int_{0}^{K}f(u)f(v)dvdu\leq
\exp(-\frac{\pi^{2}t}{8K^{2}}) .\label{upperA3}
\end{eqnarray}
Let $C_{1}=C_{1}(z)=2 (4/\pi)A(z)^{2}(1+\epsilon)^{3}$. For the
integral over $A_{2}$ we get,
\begin{eqnarray}
\ & \ & \int\!\int_{A_{2}} P_{\frac{u}{u+v}} [\eta _{(0,1)}>\frac{t}{(u+v)^{2}}] f(u)f(v)dvdu\nonumber\\
\ &\ \leq &  C_{1}\int_{0}^{K}\!\int_{K}^{\sqrt{t/M}-u} f(u)\exp (-\frac{\pi ^{2} t}{2 (u+v)^{2}}-\lambda_{D}v)dvdu\nonumber\\
 & + &  C_{1}\int_{K}^{1/2\sqrt{t/M}}\!\int_{u}^{\sqrt{t/M}-u}\sin\left( \frac{\pi u}{u+v}\right) \exp (-\frac{\pi ^{2} t}{2 (u+v)^{2}}-\lambda_{D}(u+v))dvdu \nonumber\\
\ &\ =& I+II. \label{upperA2}
\end{eqnarray}
Changing variables $u+v=z$, $u=w$
\begin{eqnarray}
\ &\ I &  = \int_{0}^{K}\!\int_{K}^{\sqrt{t/M}-u} \exp (-\frac{\pi
^{2} t}
{2 (u+v)^{2}}) f(u)\exp (-\lambda_{D}v)dvdu\nonumber\\
\ & \leq & \int_{0}^{K}\!\int_{w+K}^{\sqrt{t/M}} \exp (-\frac{\pi
^{2} t}{2 z^{2}}) f(w)\exp (-\lambda_{D}z)
\exp(\lambda_{D}w)dzdw\nonumber\\
\ & \leq & \exp(\lambda_{D}K)\int_{0}^{K}f(w)dw \int_{0}^{\infty}
\exp (-\frac{\pi ^{2} t}{2 z^{2}})
\exp (-\lambda_{D}z)dz\nonumber\\
\ & \lesssim \ & t^{1/6} \exp
(-\frac{3}{2}\pi^{2/3}\lambda_{D}^{2/3}t^{1/3}).\label{upperA21}
\end{eqnarray}
Equation (\ref{upperA21}) follows from equation
(\ref{newasymptotic}), with $a=\pi^{2}/2$, $b=\lambda_{D}$.

Changing variables $u+v=z$, $u=w$

\begin{eqnarray}
\ &\ II & \leq C_{1}\int_{K}^{1/2\sqrt{t/M}}\!\int_{2w}^{\sqrt{t/M}}
\sin\left( \frac{\pi w}{z}\right)\exp (-\frac{\pi ^{2} t}{2 z^{2}}-\lambda_{D}z)dzdw \nonumber \\
\ & \ & =C_{1}\int_{2K}^{\sqrt{t/M}}\!\int_{K}^{z/2} \sin \left(
\frac{\pi w}{z}\right)\exp (-\frac{\pi
^{2} t}{2 z^{2}}-\lambda_{D}z)dwdz \label{changevar}\\
\ & \ & \ \leq C_{1}/\pi \int_{2K}^{\sqrt{t/M}}z \cos \left(
\frac{\pi K}{z}\right) \exp (-\frac{\pi
^{2} t}{2 z^{2}}-\lambda_{D}z)dz\nonumber\\
 \ & \ & \ \leq
  (1+\epsilon)(C_{1}/\pi) 2\sqrt{\frac{\pi}{3}}(\frac{\pi^{2}}{2})^{1/2}\lambda_{D}^{-1}t^{1/2}\exp(-\frac{3}{2}\pi^{2/3}\lambda_{D}^{2/3}t^{1/3}).
\label{upperA22}
\end{eqnarray}

Equation (\ref{changevar}) follows by changing the order of the
integration. Equation (\ref{upperA22}) follows from equation
(\ref{generalast111}).

Now from equations (\ref{upperA1}), (\ref{upperA3}),
(\ref{upperA21}) and (\ref{upperA22}) we obtain
\begin{equation}\label{exactupperbound}
\limsup_{t\to\infty}t^{-1/2}\exp(\frac{3}{2}\pi^{2/3}\lambda_{D}^{2/3}t^{1/3})P_{z}[\tau
_{D}(Z)>t]\leq (1+\epsilon)(\frac{C_{1}}{\pi})
2\sqrt{\frac{\pi}{3}}(\frac{\pi^{2}}{2})^{1/2}\lambda_{D}^{-1}.
\end{equation}
Finally, from equations (\ref{exactlowerbound1}) and
(\ref{exactupperbound}) and letting $\epsilon\to 0$,
\begin{eqnarray}
 C(z)  & \leq &
\liminf_{t\to\infty} t^{-1/2}
\exp(\frac{3}{2}\pi^{2/3}\lambda_{D}^{2/3}t^{1/3})P_{z}[\tau_{D}(Z)>t]\nonumber\\
& \leq  & \limsup_{t\to\infty}
t^{-1/2}\exp(\frac{3}{2}\pi^{2/3}\lambda_{D}^{2/3}t^{1/3})
P_{z}[\tau_{D}(Z)>t]\leq   C(z),\nonumber
\end{eqnarray}
where $ C(z)=(\lambda_{D}2^{7/2})/\sqrt{3 \pi}\left(
\psi(z)\int_{D}\psi(y)dy\right) ^{2}. $
\end{proof}

\section{ Brownian-time Brownian motion in unbounded domains}\label{BTBM}
In this section we study Brownian-time Brownian motion (BTBM),
$Z_{t}^{1}$ started at $z \in \RR{R}$, in several unbounded
domains.

 If $D\subset \RR{R} ^{n}$  is an open set, write
$$
\tau_{D}(z)=\inf \{ t\geq 0: \ X_{t}+z\notin D\},
$$
and if $I\subset \RR{R}$ is an open interval, we write
$$
\eta_{I}=\inf \{ t\geq 0: \ Y_{t}\notin I\}.
$$
Let $\tau_{D}(Z^{1})$ stand for the first exit time of BTBM from $D$. We have by the continuity of paths
\begin{equation}\label{probeqisop}
P_{z}[\tau _{D}(Z^{1}) > t]=P[\eta(-\tau _{D}(z) , \tau _{D}(z) ) >t].
\end{equation}

\begin{proof}[\textbf{Proof of Theorem \ref{twisted1}}]
Let $\epsilon >0$. From Lemma \ref{isoplowasymptotic}, choose
$M>0$ so large that
$$
(1-\epsilon )\exp(-\frac{\pi ^{2} t}{8 u^{2}}) \frac{\pi t}{u^{3}}
\leq \frac{d }{du} P_{0}[\eta _{(-u,u)} > t]\leq (1+\epsilon) \
\exp(-\frac{\pi ^{2} t}{8 u^{2}}) \frac{\pi t}{u^{3}}
$$
for all $u\leq \sqrt{t/ M}$.

Let $C=C(\gamma)$. From the hypothesis choose $K>0$ so large that
\begin{eqnarray}
u^{-C(1+\epsilon)}\leq P(\tau_{D}(z) >u) \leq u^{-C(1-\epsilon) }
\ \ \mathrm{for}\ \ u\geq K. \label{dfapprox1}
\end{eqnarray}
We further assume that $t$ is so large that $K<\sqrt{t/M}$.
\begin{eqnarray}
P_{z}[\tau _{D}(Z^{1}) > t]&=&  \int_{0}^{\infty}  (\frac{d}{du}P_{0}[\eta _{(-u,u)}>t] )P[\tau _{D}(z)>u]du\nonumber\\
& \gtrsim & t \int_{K}^{\sqrt{t/M}}\exp (-\frac{\pi
^{2}t}{8u^{2}})u^{-(C(1+\epsilon)+3)}du
\end{eqnarray}
 changing variables $u^{-2}=x$, $du=-1/2x^{-3/2}dx$ the integral
 is
\begin{eqnarray}
& \gtrsim & t \int_{M/t}^{K^{-2}}\exp (-\frac{\pi ^{2}t
x}{8})x^{C(1+\epsilon)/2}dx
\end{eqnarray}
Changing variables, $z=\pi^{2}t x/8$, the integral is
\begin{eqnarray}
& \gtrsim & t^{-C(1+\epsilon)/2 } \int_{\pi^{2}
M/8}^{K^{-2}\pi^{2}t/8} e^{-z}z^{C(1+\epsilon)/2}dz.
\end{eqnarray}

Now since for some $c_{0}>0$,
$$
\int_{K^{-2}\pi^{2}t/8}^{\infty}e^{-z}z^{C(1+\epsilon)/2}dz\leq
e^{-c_{0}t},
$$
$$
\int_{0}^{\pi^{2} M/8} e^{-z}z^{C(1+\epsilon)/2}dz < \infty,
$$
and
$$
\int_{0}^{\infty} e^{-z}z^{C(1+\epsilon)/2}dz=\Gamma
(1+C(1+\epsilon)/2).
$$
 We have
\begin{eqnarray}
P_{z}[\tau _{D}(Z^{1}) > t] & \gtrsim & t^{-C(1+\epsilon)/2
}.\label{low-bound}
\end{eqnarray}
We now give an upper bound.
\begin{eqnarray}
P_{z}[\tau _{D}(Z^{1}) > t] & = & \int_{0}^{\infty}   P_{0} (\eta _{(-u,u)}>t) f(u)du\nonumber\\
& \lesssim &\int_{0}^{\sqrt{t/M}}
e^{-\frac{\pi^{2}t}{8u^{2}}}f(u)du
 + \int_{\sqrt{t/M}}^{\infty}
f(u)du\nonumber\\
& \lesssim & E\left[\exp
\left(-\frac{\pi^{2}t}{8(\tau_{D}(z))^{2}}\right)\right]
 + (\sqrt{t/M})^{-C(1-\epsilon)}\nonumber\\
 & \lesssim &  t^{-C(1-\epsilon)/2 } \label{upp-bound}
\end{eqnarray}
Equation (\ref{upp-bound}) follows from Lemma \ref{upper-bound}
and the asymptotics of $\tau_{D}(z)$.

\noindent Now from Equations (\ref{low-bound}) and
(\ref{upp-bound}) we have

\begin{eqnarray}
t^{-C(1+\epsilon)/2} \lesssim P_{z}[\tau _{D}(Z^{1}) > t] &
\lesssim & t^{-C(1-\epsilon)/2 }.\nonumber
\end{eqnarray}
Now taking logarithm of the above inequalities and then dividing
by $\log t$ and letting $\epsilon\to 0$, we obtain the desired
result.
\end{proof}

\begin{proof}[\textbf{Proof of Theorem \ref{theorem3}}]
Let $\epsilon >0$. From Lemma \ref{isoplowasymptotic}, choose
$M>0$ so large that
\begin{equation}\label{du-bound}
(1-\epsilon )\exp(-\frac{\pi ^{2} t}{8 u^{2}}) \frac{\pi t}{u^{3}}
\leq \frac{d }{du} P_{0}[\eta _{(-u,u)} > t]\leq (1+\epsilon) \
\exp(-\frac{\pi ^{2} t}{8 u^{2}}) \frac{\pi t}{u^{3}}
\end{equation}
for all $u\leq \sqrt{t/ M}$.

Let $C=j^{2}_{\nu}$. From the hypothesis choose $K>0$ so large that
\begin{eqnarray}
& & \exp(-C(1+\epsilon)u(\log u)^{-2/p} )\nonumber\\
& & \leq P(\tau_{P_{f}}(z) >u) \nonumber\\
 & & \leq \exp(-C(1-\epsilon)u(\log u)^{-2/p} )
\label{dfapprox2}
\end{eqnarray}
for $ u\geq K.$
We further assume that $t$ is so large that $K<\sqrt{t/M}$.

Then, by equations (\ref{du-bound}) and (\ref{dfapprox2})

$ P_{z}[\tau _{D}(Z^{1}) > t]$
\begin{eqnarray}
 & \gtrsim & t
\int_{K}^{\sqrt{t/M}}u^{-3}\exp (-\frac{\pi
^{2}t}{8u^{2}})\exp(-C(1+\epsilon)u(\log u)^{-2/p} )du
\end{eqnarray}
 changing variables $u^{-2}=x$, $du=-1/2x^{-3/2}dx$ the integral
 is
\begin{eqnarray}
& \gtrsim & t \int_{M/t}^{K^{-2}}\exp (-\frac{\pi ^{2}t
x}{8})\exp(-C(1+\epsilon)x^{-1/2}(\log x^{-1/2})^{-2/p} )dx\nonumber
\end{eqnarray}
 Now we set
 $$
 v(x)=C(1+\epsilon)x^{-1/2}(\log x^{-1/2})^{-2/p}
 $$
which gives
$$
dv=C(1+\epsilon)x^{-3/2}(\log x^{-1/2})^{-2/p}[-1/2+1/p(\log
x^{-1/2})^{-1}]dx= -Vdx
$$

Then the integral is

\begin{eqnarray}
& \gtrsim & t \int_{M/t}^{K^{-2}}VV^{-1}\exp (-\frac{\pi ^{2}t
x}{8})\exp(-C(1+\epsilon)x^{-1/2}(\log x^{-1/2})^{-2/p} )dx\nonumber
\end{eqnarray}

Now
$$
V^{-1}\gtrsim t^{-3/2}[1/2-1/p(\log \sqrt{t/M})^{-1}]^{-1}\gtrsim
t^{-1/2}
$$
Hence the integral is
\begin{eqnarray}
& \gtrsim & t^{1/2} \int_{M/t}^{K^{-2}}\exp (-\frac{\pi ^{2}t
x}{8})\exp(-v )Vdx
\end{eqnarray}
Now from Lemma \ref{densitylemma} and de Bruijn's Tauberian Theorem \ref{debruijn} we have
$$
\log \int_{0}^{\infty}\exp (-\frac{\pi ^{2}t x}{8})\exp(-v )(-dv)
\sim -C(p,\epsilon)t^{1/3}(\log t)^{-4/(3p)}
$$
where $C(p,\epsilon)=(3/2)^{(4+3p)/3p}2^{1/3}((1+\epsilon)j_{\nu}^{2}2^{2/p})^{2/3}(\pi ^{2}/8)^{1/3}$.

From the bounds for some $c_{1}, c_{2}>0$,
$$
\int_{0}^{M/t}\exp (-\frac{\pi ^{2}t x}{8})\exp(-v )(-dv) \leq
e^{-c_{1}t^{1/2} (\log t)^{-2/p} },
$$
and
$$
\int_{K^{-2}}^{\infty}\exp (-\frac{\pi ^{2}t x}{8})\exp(-v )(-dv)
\leq e^{-c_{2}t },
$$
 we get
\begin{equation}\label{lower-f}
P_{z}[\tau _{P_{f}}(Z^{1}) > t] \gtrsim \exp(-(1+\epsilon)^{2}C(p,\epsilon)t^{1/3}(\log t)^{-4/3p})
\end{equation}
We next give the upper bound
\begin{eqnarray}
P_{z}[\tau _{P_{f}}(Z^{1}) > t] & = & \int_{0}^{\infty}   P_{0} (\eta _{(-u,u)}>t) f(u)du\nonumber\\
& \lesssim &\int_{0}^{\sqrt{t/M}}
e^{-\frac{\pi^{2}t}{8u^{2}}}f(u)du
 + \int_{\sqrt{t/M}}^{\infty}
f(u)du\nonumber\\
& \lesssim & E\left[\exp
\left(-\frac{\pi^{2}t}{8(\tau_{P_{f}}(z))^{2}}\right)\right]
 + \exp(-C(1-\epsilon)t^{1/2}(\log t)^{-2/p})\nonumber
\end{eqnarray}
The upper bound follows from de Bruijn's Tauberian Theorem by observing that  as $x\to 0^{+}$
$$
\log P[1/(\tau_{P_{f}}(z))^{2}\leq x]= \log P[\tau_{P_{f}}(z)\geq
x^{-1/2}]\sim -Cx^{-1/2}(|\log x^{-1/2}|)^{-2/p}
$$
Hence
\begin{equation}\label{upper-f}
P_{z}[\tau _{P_{f}}(Z^{1}) > t] \lesssim \exp(-(1-\epsilon)D(p,\epsilon)t^{1/3}(\log t)^{-4/3p})
\end{equation}
where $D(p,\epsilon)=(3/2)^{(4+3p)/3p}2^{1/3}((1-\epsilon)j_{\nu}^{2}2^{2/p})^{2/3}(\pi ^{2}/8)^{1/3}.$

Therefore from equations  (\ref{lower-f}) and (\ref{upper-f}), we obtain
\begin{eqnarray}
& & \exp(-(1+\epsilon)^{2}C(p,\epsilon)t^{1/3}(\log t)^{-4/3p})\nonumber\\
& & \lesssim P_{z}[\tau _{P_{f}}(Z^{1}) > t] \nonumber\\
& & \lesssim \exp(-(1-\epsilon)^{2}D(p,\epsilon)t^{1/3}(\log t)^{-4/3p})\nonumber
\end{eqnarray}
Now taking logarithms and letting $\epsilon\to 0$, we get the desired result.
\end{proof}

\begin{proof}[\textbf{Proof of Theorem \ref{not-sharp}}]
The proof follows the same steps of the proof of Theorem \ref{theorem3}, so we omit the details.
 \end{proof}

\section{Asymptotics}\label{approximation}

In this Section we will recall some lemmas that were used in
section 3 and section 4. The following lemma is proved in
\cite[Lemma A1]{deblassie} (it also follows from more general
results on
``intrinsic ultracontractivity").  We include it for completeness.

\begin{lemma}\label{lemmaA.1}
As $t\rightarrow \infty$,
\[
P_{x}[\eta _{(0,1)} >t]\  \sim \  \frac{4}{\pi}e^{-\frac{\pi ^{2}
t}{2} } \sin \pi x, \ \  \mathrm{uniformly} \ \ \mathrm{for}\ \
x\in (0,1).
\]
\end{lemma}

We  recall a result from Nane \cite[Lemma
6.2]{nane} that will be used for the process $Z^{1}$.
\begin{lemma}\label{isoplowasymptotic} Let $B=\{ u>0: \ t/ u^{2} >M\}$ for $M$ large.
 Then on $B$,
\begin{equation}
 \frac{d }{du} P_{0}[\eta _{(-u,u)} > t]\sim \
\exp(-\frac{\pi ^{2} t}{8 u^{2}}) \frac{\pi t}{u^{3}}.
\label{lowerdiffisop}
\end{equation}
\end{lemma}

\textbf{Acknowledgments.} I would like to thank  Professor Rodrigo
 Ba\~{n}uelos, my academic advisor, for suggesting this problem to me and for  his guidance on this paper.

\end{document}